\input amstex
\documentstyle{amsppt}
\magnification=\magstep1
\NoBlackBoxes

\def\dim{\mathop{\roman{dim}}}
\def\int{\mathop{\roman{int}}}

\def\1{^{-1}}

\redefine\int{\operatorname{Int}}
\define\ExD{\operatorname{ext--dim}}

\redefine\int{\operatorname{Int}}

\hsize=6.5truein
\topmatter
\title
Extension dimension for paracompact spaces
\endtitle
\author
Jerzy Dydak
\endauthor
\address University of Tennessee, Knoxville, TN 37996
\endaddress
\email
dydak$\@$math.utk.edu
\endemail
\date 
October 26, 2002
\enddate
\thanks Research supported in part by a grant
 DMS-0072356 from the National Science Foundation
\endthanks
\abstract
We prove existence of extension dimension for paracompact spaces.
Here is the main result of the paper:
\proclaim{Theorem} Suppose $X$ is a paracompact space.
There is a CW complex $K$ such that
\item{a.} $K$ is an absolute extensor of $X$ up to homotopy,
\item{b.} If a CW complex $L$ is an absolute extensor of
$X$ up to homotopy, then $L$ is an absolute extensor of
$Y$ up to homotopy of any paracompact space $Y$ such that
$K$ is an absolute extensor of
$Y$ up to homotopy.
\endproclaim
The proof is based on the following simple result (see 1.6).
\proclaim{Theorem} Suppose $X$ be a paracompact space
and $f:A\to Y$ is a map from a closed subset $A$ of $X$
to a space $Y$. $f$ extends over $X$ if $Y$ is
the union of a family $\{Y_s\}_{s\in S}$ of its subspaces with the following
properties:
\item{a.} Each $Y_s$ is an absolute extensor of $X$,
\item{b.} For any two elements $s$ and $t$ of $S$ there is $u\in S$
such that $Y_s\cup Y_t\subset Y_u$,
\item{c.} 
$A=\bigcup\limits_{s\in S} \int_A(f^{-1}(Y_s))$.
\endproclaim
That result implies a few well-known theorems of classical theory of retracts
which makes it of interest in its own.
\endabstract
\dedicatory  
Dedicated to Jed Keesling on the occasion of his sixtieth birthday.
 \enddedicatory
\keywords dimension, extension of maps, absolute extensors, CW complexes, paracompact spaces
\endkeywords
\subjclass 54C55, 54F45
\endsubjclass
\endtopmatter

\document

\head 1. Introduction \endhead
A. Dranishnikov \cite{Dr} introduced the concept of extension dimension
for compact Hausdorff spaces as a generalization of both covering dimension
and cohomological dimension.
\definition{1.1. Definition} Suppose $X$ is a compact Hausdorff space.
A CW complex $K$ is called the {\bf extension dimension} of $X$ 
if the following two conditions are satisfied:
\item{a.} $K$ is an absolute extensor of $X$,
\item{b.} If a CW complex $L$ is an absolute extensor of
$X$, then $L$ is an absolute extensor of
$Y$ for any compact Hausdorff space $Y$ such that
$K$ is an absolute extensor of
$Y$.
\enddefinition

The meaning of Definition 1.1 is that extension dimension
of $X$ is the minimal element of a subclass in a certain order on the class
of all CW complexes. Namely, one can define $K\leq L$
if $\Cal C_K\subset \Cal C_L$, where $\Cal C_M$ is the class
of all compact Hausdorff spaces $X$ such that $M\in AE(X)$.
Now, $K$ is the extension dimension of $X$ if it is the minimal
element among all $L$ such that $X\in \Cal C_L$.

\par One can ponder the existence of extension dimension for
other classes of topological spaces. This was done 
by A.Dranishnikov and J.Dydak in \cite{D-D$_1$} for separable
metrizable spaces, and by I.Ivan\v si\' c and  
L.Rubin in \cite{I-R} for metrizable spaces. However,
the proofs in \cite{D-D$_1$} and \cite{I-R} are quite
complicated. The author believes that, for a theory to be successful, its
foundations should be fairly simple. The purpose of this paper
is to provide quite an elementary proof of the existence
of extension dimension for paracompact spaces.

\par One of the main ideas of extension theory is to
investigate spaces by mapping them (or their subspaces) to spaces $K$ with good
local properties. Traditionally, the spaces one wants to investigate
are metrizable or compact Hausdorff. That tradition is the result of 
a natural evolution:
euclidean spaces, their subspaces, their compactifications.
Also, two classes of spaces with good local properties
emerged; CW complexes and ANRs (absolute neighborhood retracts
of metrizable spaces). Those two classes are known to be identical
up to homotopy but as of now we do not know of a single
class which could be used in their place.
Is there a natural class of spaces which naturally combines metrizable
spaces and compact Hausdorff spaces?
The problem is that ANRs do not have to be absolute
neighborhood extensors of compact Hausdorff spaces.
One could bypass that problem by
considering only maps $f:A\to K$ on
closed subsets $A$ of $X$ which are $G_\delta$-subsets
of $X$. Since being closed and a $G_\delta$ subset
of a normal space $X$ is equivalent to be a zero subset
(i.e., a set of the form $\alpha^{-1}(0)$ for some
continuous $\alpha:X\to [0,1]$), let us formulate
the corresponding variation of the concept
of absolute extensor.

\definition{1.2. Definition} 
$Y\in AE_0(X)$ ($Y\in ANE_0(X)$, respectively)
means that all maps $f:A\to Y$
extend over $X$ (over a neighborhood of $A$ in $X$, respectively)
provided $A$ is a zero subset of $X$.
\enddefinition

It is known that, if $K$ is an ANR
and $X$ is paracompact space, then $K\in ANE_0(X)$.
However, if $K$ is a CW complex the analogous statement is false.
Indeed, van Douwen and Pol \cite{D-P} constructed
the strongest possible counterexample. In their case (see section 3
of \cite{D-P})
$K$ is the cone over infinite discrete CW complex
and $A$ is a closed subspace of a countable paracompact space $X$.
\par To avoid problems with extending maps to CW complexes
over neighborhoods of closed subsets of paracompact spaces the papers \cite{D-D$_1$} and \cite{I-R}
 create subclasses of paracompact spaces.
In \cite{D-D$_1$} cw-spaces are defined as paracompact $k$-spaces $X$
such that any contractible CW complex $K$ is an absolute extensor of $X$.
In \cite{I-R} dd-spaces are defined. 
\par In this paper the difficulty is avoided by
switching the focus from extending maps to extending maps
up to homotopy which may seem to be a more difficult
task. However, there is a special class of generic maps
to CW complexes (called locally compact maps)
for which the two extension problems are equivalent.
As a result we obtain three possible interpretations
of extension dimension for paracompact spaces:

\proclaim{1.3. Theorem} Suppose $X$ is a paracompact space.
There is a CW complex $K$ (called the {\bf extension dimension} of $X$)
such that
\item{a.} $K$ is an absolute extensor of $X$ up to homotopy,
\item{b.} If a CW complex $L$ is an absolute extensor of
$X$ up to homotopy, then $L$ is an absolute extensor of
$Y$ up to homotopy of any paracompact space $Y$ such that
$K$ is an absolute extensor of
$Y$ up to homotopy.
\endproclaim

\proclaim{1.4. Theorem} Suppose $X$ is a paracompact space.
There is a simplicial complex $K$ such that
\item{a.} $|K|_m$ is an absolute extensor of $X$ and is complete,
\item{b.} If a complete ANR $L$ is an absolute extensor of
$X$, then $L$ is an absolute extensor 
any paracompact space $Y$ such that
$|K|_m$ is an absolute extensor of
$Y$.
\endproclaim

\proclaim{1.5. Theorem} Suppose $X$ is a paracompact space.
There is a simplicial complex $K$ such that
\item{a.} $|K|_m\in AE_0(X)$,
\item{b.} If $L\in AE_0(X)$ is an ANR, then $L\in AE_0(Y)$ 
for any paracompact space $Y$ such that
$|K|_m\in AE_0(Y)$.
\endproclaim

\par
Let us start with
a general, yet simple, result which is at the core of our approach
to extension dimension theory.

\proclaim{1.6. Theorem} Suppose $X$ be a paracompact space
and $f:A\to Y$ is a map from a closed subset $A$ of $X$
to a space $Y$. $f$ extends over $X$ if $Y$ is
the union of a family $\{Y_s\}_{s\in S}$ of its subspaces with the following
properties:
\item{a.} Each $Y_s$ is an absolute extensor of $X$,
\item{b.} For any two elements $s$ and $t$ of $S$ there is $u\in S$
such that $Y_s\cup Y_t\subset Y_u$,
\item{c.} 
$A=\bigcup\limits_{s\in S} \int_A(f^{-1}(Y_s))$.
\endproclaim
\demo{Proof} Define $U_s=(X-A)\cup \int_A(f^{-1}(Y_s))$ for each $s\in S$. 
Each $U_s$ is an open subset of $X$
and $X=\bigcup\limits_{s\in S}U_s$. Since $X$ is paracompact, there is
a locally finite partition of unity $\{g_s\}_{s\in S}$ on $X$
such that $g_s^{-1}(0,1]\subset U_s$ for each $s\in S$ (see \cite{En}, Theorem 5.1.9
and its proof).
For all finite subsets $T$ of $S$ define
$B_T=\{x\in X\mid g_s(x)>0\implies s\in T\}$.
We plan to create, for all finite subsets $T$ of $S$, 
elements $a(T)$ of $S$ and maps $f_T:B_T\to Y_{a(T)}$ so that
the following conditions are satisfied:
\item{1.} $Y_{a(F)}\subset Y_{a(T)}$ for each $F\subset T$,
\item{2.} $f_T|B_F=f_F$ for each $F\subset T$,
\item{3.} $f_T|A\cap B_T=f|A\cap B_T$.
\par This is going to be accomplished by induction on the number of elements of $T$.
For one-element sets $T=\{s\}$ we simplify notation to $T=s$.
Notice that $B_s=g_s^{-1}(1)$ for each $s\in S$. $\{B_s\}_{s\in S}$ is a discrete family
and $f(A\cap B_s)\subset Y_s$ for each $s\in S$. Therefore we can extend
each $f|A\cap B_s$ to $f_s:B_s\to Y_s$ and we put $a(s)=s$.
Suppose $f_T$ and $a(T)$ exist for all $T$ with cardinality at most $n$.
Given $T$ containing exactly $n+1$ elements, pick $s\in S$ so that
$Y_s$ contains all of $Y_{a(F)}$ with $F$ being a proper subset of $T$. Put $a(T)=s$.
All of $f_F$, $F$ a proper subset of $T$, can be pasted together
and produce a map $h$ on a closed subset $B$ of $B_T$ with values in $Y_s$
and extending $f$ on $A\cap B$. Since $f(A\cap B_T)\subset Y_s$,
$h$ extends over $B_T$ and produces $f_T:B_T\to Y_{a(T)}$ with the
desired properties.
\par Since $B_T\cap B_F=B_{T\cap F}$, all $f_T$ can be pasted together
to produce a function $f':X\to Y$ which is an extension of $f$.
Any point $x\in X$ has a neighborhood $U$ which intersects
only finitely many of $g_s^{-1}(0,1]$ which means that there is a finite set
 $T$ such that $U\subset  B_{T}$.
As $f'|B_{T}$ is continuous, so is $f'|U$
which completes the proof.
\qed
\enddemo

Before applying 1.6 let us recall a canonical method from \cite{Dy$_2$}
of converting results about absolute extensors
to theorems about absolute neighborhood extensors.
This is done by using the so-called covariant cones.
For any space $P$ its {\bf covariant cone}
 $Cone(P)$ is $P\times I/P\times \{1\}$
with the topology induced by open sets in $P\times [0,1)$
and a basis of neighborhoods of the vertex $v=P\times \{1\}/P\times \{1\}$
being $P\times (t,1]/P\times \{1\}$, $t\in [0,1)$.
 In \cite{Dy$_2$} (see Theorem 2.9)
it is shown that if $P$ is Hausdorff, contains at least two points, and is an absolute
neighborhood extensor of a space $M$, then $Cone(P)$ is an absolute
extensor of $M$. Notice that, in case of normal spaces $M$, the proof
of 2.9 in \cite{Dy$_2$} applies to all spaces $P$ as the assumption
of $P$ being Hausdorff and containing at least two points was used only to deduce
that $M$ is normal.

\proclaim{1.7. Corollary} Suppose $X$ be a paracompact space
and $f:A\to Y$ is a map from a closed subset $A$ of $X$
to a space $Y$. $f$ extends over a neighborhood of $A$ in $X$ if $Y$ is
the union of a family $\{Y_s\}_{s\in S}$ of its subspaces with the following
properties:
\item{a.} Each $Y_s$ is an absolute neighborhood extensor of $X$,
\item{b.} For any two elements $s$ and $t$ of $S$ there is $u\in S$
such that $Y_s\cup Y_t\subset Y_u$,
\item{c.} 
$A=\bigcup\limits_{s\in S} \int_A(f^{-1}(Y_s))$.
\endproclaim
\demo{Proof}
Let $Z=Cone(Y)$ with vertex $v$ and $Z_s=Cone(Y_s)$ for each $s\in S$. 
 Therefore, $f$ considered as a map from $A$
to $Z$ satisfies hypotheses of Theorem 1.6 and extends over $X$. Let $g:X\to Z$ be an extension
of $f$
and let $U=g^{-1}(Z-\{v\})$. There is a retraction $r:Z-\{v\}\to Y$
which means that the composition of $g|U$ and $r$ produces an extension
$f':U\to Y$ of $f$.
\qed
\enddemo

The strength of 1.7 is that it implies two well-known results
from the theory of retracts and its proof is much simpler than those of
original results. The first one is a theorem first proved by Dugundji \cite{Du} 
(and independently by Kodama \cite{Ko})
 for the special case of simplicial complexes with the CW topology.
In full generality it follows from a result of Cauty \cite{Ca}
that each CW complexes $K$ can be embedded in a polyhedron with CW topology
in such a way that there is a retraction $r:U\to K$ from a neighborhood
$U$ of $K$.
\proclaim{1.8. Corollary (Cauty-Dugundji-Kodama)} CW complexes are absolute neighborhood extensors
of metrizable spaces.
\endproclaim
\demo{Proof} Finite subcomplexes of a CW complex $K$
form a family closed under finite sums, each of them is an absolute neighborhood
extensor
of normal spaces, and any map $f:A\to K$ from a first countable space
has the property that each point $x\in A$ has a neighborhood $U$
such that $f(U)$ is contained in a finite subcomplex of $K$
(see \cite{Dy$_2$}, Corollary 4.5). Thus, 1.7 applies.
\qed
\enddemo

The second one is a result of Hanner as proved in \cite{Hu} in quite a complicated way
 on eleven pages (see Theorem 17.1 on pp. 68--79).
\proclaim{1.9. Corollary (Hanner)} Suppose $X$ is a paracompact space. If a Hausdorff space
$Y$
is a union of open subsets $U$ which are absolute neighborhood extensors of
$X$, then $Y$ is an absolute neighborhood extensor of $X$.
\endproclaim
\demo{Proof} The family of all open subsets of $Y$ which are absolute neighborhood
extensors of $X$ is closed under finite unions (see \cite{Hu}, Theorem 8.2),
so 1.7 applies.
\qed
\enddemo

The author would like to thank Sergey Antonyan for asking questions about existence of a simple
proof of Cauty-Dugundji-Kodama Theorem 1.8, and to Ivan Ivan\v si\' c for help with sorting out
the issues related to CW complexes and ANE for paracompact spaces.
Antonyan's question stemmed from \cite{AEM}, where a proof of 1.8
is given which is simpler than the original one. Also, it is mentioned in \cite{AEM}
that our approach, when applied to the equivariant case, is of interest
and offers simplifications similar to those in the non-equivariant case.

\head 2. Locally compact maps \endhead

The simplicity of 1.6-1.7 and their applications made the author
think that one should attempt to build extension theory based on 1.6.
Since our interest is mostly in maps to CW complexes, the proof of 1.8
suggests that we need to concentrate on maps such that
every point has a neighborhood whose image is contained in a finite
subcomplex. A generalization to arbitrary spaces is obvious:
\definition{2.1. Definition} A map $f:X\to Y$ is called {\bf locally compact}
if for every element $x\in X$ there is a neighborhood $U$ in $X$
such that $f(U)$ is contained in a compact subset of $Y$.
\enddefinition

\remark{Remark} It is easy to show that $f:X\to Y$ is locally compact
if and only if for any compact subset $Z$ of $X$ there is a neighborhood
$U$ of $Z$ in $X$ such that $f(U)$ is contained in a compact subset of $Y$.
\endremark
Let us point out that, in the case of maps to simplicial complexes
with the weak topology, the concept of locally compact map
corresponds to the concept of locally finite partition of unity.
In 2.2 and in the remainder of the paper we follow the notation of \cite{M-S},
where $|L|_w$ is the body of a simplicial complex $L$ equipped with the weak topology,
and $|L|_m$ is the body of a simplicial complex $L$ equipped with the metric topology.

\proclaim{2.2. Proposition} Let $L$ be a simplicial complex.
A map $f:X\to |L|_w$ is locally compact if and only if
the corresponding partition of unity on $X$ is locally finite.
\endproclaim
\demo{Proof} Let $V$ be the set of vertices of $L$. The partition of unity
corresponding to $f$ is the set of maps $f_v:X\to I$ (those are the barycentric coordinates
of $f(x)$ according to the terminology of \cite{M-S})
so that $f(x)=\sum\limits_{v\in V}f_v(x)\cdot v$. $\{f_v\}_{v\in V}$ being locally finite
means that each point $x\in X$ has a neighborhood $U$ such that only finitely many
$f_v$ are non-zero on $U$. That is the same as saying that $f(U)$ is contained
in a finite subcomplex of $|L|_w$.
\qed
\enddemo

The remainder of this section is devoted to the homotopy theory of locally compact maps.
We start with a few elementary observations.
\proclaim{2.3. Proposition} Suppose $f:X\to Y$ and $g:Y\to Z$ are maps.
If $f$ or $g$ is locally compact, then $g\circ f$ is locally compact.
\endproclaim
\demo{Proof} Suppose $x\in X$. If there is a neighborhood $U$
of $x$ in $X$ such that $f(U)$ is contained in a compact
subset $C$ of $Y$, then $gf(U)$ is contained in $g(C)$ which is compact.
If $f(x)$ is contained in a neighborhood $V$ in $Y$ such that
$g(V)$ is contained in a compact subset $C$ of $Z$,
then we put $U=g^{-1}(V)$ and notice that $gf(U)$ is contained in $C$.
\qed
\enddemo

\proclaim{2.4. Proposition} Suppose $X$ is the union of
a locally finite family $\{X_s\}_{s\in S}$ consisting of closed sets.
Let $f:X\to Y$
be a map. If $f|X_s$ is locally compact for each $s\in S$, 
then $f$ is locally compact.
\endproclaim
\demo{Proof} Suppose $x\in X$. If $x\in X_s$ for some $s\in S$, we pick
a neighborhood $U_s$ of $x$ in $X$ such that $f(U_s\cap X_s)$
is contained in a compact subset $C_s$ of $Y$.
Let $T$ be a finite subset of $S$
such that $x\in X_s$ if and only if $s\in T$.
Let $W=X-\bigcup\limits_{s\in S-T} X_s$, and put
$U=W\cap \bigcap\limits_{s\in T} U_s$.
Obviously, $U$ is a neighborhood of $x$ in $X$.
It remains to show that $f(U)\subset \bigcup\limits_{s\in T} C_s$
which follows
from $U\subset \bigcup\limits_{s\in T} U_s\cap X_s$.
\qed
\enddemo

\proclaim{2.5. Proposition} If $f_i:X_i\to Y_i$ is locally compact for $i=1,2$, 
then $f_1\times f_2:X_1\times X_2\to Y_1\times Y_2$ is locally compact.
\endproclaim
\demo{Proof} Suppose $(x_1,x_2)\in X_1\times X_2$.
Pick a neighborhood $U_i$ of $x_i$ in $X_i$ such that
$f_i(U_i)$ is contained in a compact subset $C_i$ of $Y_i$, $i=1,2$.
Notice that $(f_1\times f_2)(U_1\times U_2)\subset C_1\times C_2$
and $C_1\times C_2$ is compact.
\qed
\enddemo

Our next two results show that locally compact maps are prevalent,
up to homotopy, among maps to CW complexes.
\proclaim{2.6. Proposition} If $X$ is homotopy equivalent to a CW complex,
then $id_X:X\to X$ is homotopic to a locally compact map.
\endproclaim
\demo{Proof} First consider $X=|L|_w$, where $L$ is a simplicial complex.
$X$ is paracompact and open stars $\{St(v,L)\}$, $v$ is a vertex of $L$, form an open cover
of $X$. Therefore we can find a locally finite partition of unity
$\{g_v\}$ on $X$ so that $g_v^{-1}(0,1]\subset St(v,L)$ for each $v$
(see \cite{En}, Lemma 5.1.8, Theorem 5.1.9 and its proof).
That partition of unity induces a locally compact map $g:X\to X$
with the property that if $x$ belongs to a simplex $\Delta$, then $g(x)\in \Delta$.
The function $H:X\times I\to X$ defined by $H(x,t)=(1-t)\cdot x+t\cdot g(x)$
is continuous on $\Delta\times I$ for each simplex $\Delta$ which means
that $H$ is continuous. Thus, $H$ is a homotopy joining $id_X$ and $g$.
\par 	If $X$ is homotopy equivalent to a CW complex, then we can find maps
$u:X\to Y=|L|_w$ and $d:|L|_w\to X$ such that $d\circ u$ is homotopic
to the identity $id_X$ (see \cite{M-S}). Let $h:Y\to Y$ be a locally compact map
homotopic to $id_Y$. Put $g=d\circ h\circ u$. Notice that $g$
is a locally compact map (use 2.3) homotopic to $id_X$.
\qed
\enddemo

\proclaim{2.7. Corollary} Suppose $Y$ is a space such that $id_Y:Y\to Y$ is homotopic to
a locally compact map. If $f:X\to Y$ is a map such that $f|A$ is locally compact
for some closed subset $A$ of $X$, then there is a homotopy $H:X\times I\to Y$ starting at $f$
such that $H|A\times I\cup Y\times \{1\}$ is locally compact.
\endproclaim
\demo{Proof} Let $G:Y\times I\to Y$ be a homotopy joining $id_Y$ and
a locally compact map. Define $H$ as $G\circ (f\times id_I)$.
$H$ starts at $f$, $H|X\times \{1\}$ is the composition of $f$ and a locally compact map,
and $H|A\times I$ is the composition of $f\times id_I|A\times I$ (which is a locally
compact map by 2.4)
and $H|A\times I$. By 2.3 and 2.4, $H|A\times I\cup Y\times \{1\}$ is locally compact.
\qed
\enddemo

Our strategy from now on is to replace every map by a homotopic locally
compact map. That calls for obvious generalizations of well-known concepts
which will be useful in simplifying the exposition.
\definition{2.8. Definition} Suppose $X$ is a space and $K$ is a CW complex.
$K\in AE_{lc}(X)$ means that any locally compact map $f:A\to K$
on a closed subset $A$ of $X$ extends to a locally compact map 
$f':X\to K$.
\enddefinition

We are now ready for an analog of 1.6 which will be our main tool
in presenting the extension theory of paracompact spaces.
\proclaim{2.9. Theorem}  Suppose a 
CW complex $K$ is
the union of a family $\{K_s\}_{s\in S}$ of its subcomplexes so that 
for any two elements $s$ and $t$ of $S$ there is $u\in S$
with $K_s\cup K_t\subset K_u$.
 Let $X$ be a paracompact space. If, for each $s\in S$,
there is $t\in S$ so that any locally compact map $f:A\to K_s$
from a closed subset $A$ of $X$ extends to a locally compact
map $f':X\to K_t$,
then $K\in AE_{lc}(X)$.
\endproclaim
\demo{Proof} Suppose $f:A\to K$ is a locally compact map,
where $A$ is a closed subset of $X$. Given $x\in A$ there is
a neighborhood $U$ of $x$ in $A$ so that $f(U)$ is contained
in a compact subset $Z$ of $K$. Each compact subset of a CW complex
is contained in a finite subcomplex which must be contained
in $K_s$ for some $s\in S$. Therefore interiors (in $A$) of
sets $f^{-1}(K_s)$ cover $A$.
\par Define $U_s=(X-A)\cup \int_A(f^{-1}(K_s))$ for each $s\in S$. 
Each $U_s$ is an open subset of $X$
and $X=\bigcup\limits_{s\in S}U_s$. Since $X$ is paracompact, there is
a locally finite partition of unity $\{g_s\}_{s\in S}$ on $X$
such that $g_s^{-1}(0,1]\subset U_s$ for each $s\in S$ 
(see \cite{En}, Lemma 5.1.8, Theorem 5.1.9 and its proof).
For all finite subsets $T$ of $S$ define
$B_T=\{x\in X\mid g_s(x)>0\implies s\in T\}$.
We plan to create, for all finite subsets $T$ of $S$, the objects
\item{i.}
elements $a(T), b(T)$ of $S$, 
\item{ii.} locally compact maps $f_T:B_T\to K_{b(T)}$ 
\par\noindent
so that
the following conditions are satisfied:
\item{1.} $K_{a(F)}\subset K_{a(T)}$ for each $F\subset T$,
\item{2.} $K_{b(F)}\subset K_{b(T)}$ for each $F\subset T$,
\item{3.} any locally compact map $h:D\to K_{a(T)}$
on a closed subset $D$ of $X$ extends to a locally compact map $h':X\to K_{b(T)}$,
\item{4.} $f_T|B_F=f_F$ for each $F\subset T$,
\item{5.} $f_T|A\cap B_T=f|A\cap B_T$.
\par This is going to be accomplished by induction on the number of elements of $T$.
For one-element sets $T=\{s\}$ we simplify notation to $T=s$.
Notice that $B_s=g_s^{-1}(1)$ for each $s\in S$. $\{B_s\}_{s\in S}$ is a discrete family
and $f(A\cap B_s)\subset K_s$ for each $s\in S$. We put $a(s)=s$
and we find $t=b(s)$ so that
any locally compact map $h:D\to K_{s}$
on a closed subset $D$ of $X$ extends to a locally compact map $h':X\to K_{t}$.
Therefore we can extend
each $f|A\cap B_s$ to a locally compact
$f_s:B_s\to K_t$.
Suppose $f_T$, $a(T)$, and $b(T)$ exist for all $T$ with cardinality at most $n$.
Given $T$ containing exactly $n+1$ elements pick $s\in S$ so that
$K_s$ contains all of $K_{b(F)}$ with $F$ being a proper subset of $T$. Put $a(T)=s$.
We find $t=b(T)$ so that
any locally compact map $h:D\to K_{s}$
on a closed subset $D$ of $X$ extends to a locally compact map $h':X\to K_{t}$.
All of $f_F$, $F$ a proper subset of $T$, can be pasted together 
and produce a locally compact (see 2.4) map $h$ on a closed subset $B$ of $B_T$ with values in $K_s$
and extending $f$ on $A\cap B$. Since $f(A\cap B_T)\subset K_s$,
$h$ extends over $B_T$ and produces $f_T:B_T\to K_{b(T)}$ with the
desired properties.
\par Since $B_T\cap B_F=B_{T\cap F}$, all $f_T$ can be pasted together
to produce a function $f':X\to K$ which is an extension of $f$.
Any point $x\in X$ has a neighborhood $U$ which intersects
only finitely many of $g_s^{-1}(0,1]$ which means that there is a finite set
 $T$ such that $U\subset B_{T}$.
As $f'|B_{T}$ is locally compact, so is $f'|U$
which completes the proof.
\qed
\enddemo

\proclaim{2.10. Corollary} If $X$ is a paracompact space and $K$ is
a contractible CW complex, then $K\in AE_{lc}(X)$.
\endproclaim
\demo{Proof} Consider the cone $Cone(K)$ of $K$ with the
weak topology. The family of cones of finite subcomplexes of $K$
forms a family satisfying hypotheses of 2.9. Since $K$ is a retract
of its cone, $K\in AE_{lc}(X)$.
\qed
\enddemo

Our next result says that CW complexes are absolute neighborhood
extensors of paracompact spaces if the class of locally
compact maps is considered (notice that it does not make sense to talk about
category of locally compact maps as identity $id_X:X\to X$
is locally compact if and only if $X$ is locally compact).

\proclaim{2.11. Corollary} If $X$ is a paracompact space, $K$ is
a CW complex, and $f:A\to K$ is a locally compact map on a closed 
subset $A$ of $X$,
then there exists a locally compact extension $f':U\to K$ of $f$
over a neighborhood $U$ of $A$ in $X$.
\endproclaim
\demo{Proof} By 2.10 any locally compact map $f:A\to K$,
$A$ closed in $X$, extends to a locally compact $g:X\to Cone(K)$.
Let $v$ be the vertex of $Cone(K)$. Put $U=g^{-1}(Cone(K)-\{v\})$,
$r:Cone(K)-\{v\}\to K$ the canonical retraction, and $f'=r\circ (g|U)$.
\qed
\enddemo

We will also need a Homotopy Extension Theorem for locally compact maps.

\proclaim{2.12. Corollary} Suppose $X$ is a paracompact space,
$A$ is a closed subset of $X$, and $K$ is
a CW complex. If $H:A\times I\cup X\times \{0\}\to K$ is a locally
compact map, then it extends to a locally compact $H':X\times I\to K$.
\endproclaim
\demo{Proof} By 2.11 there is an open neighborhood $V$ of 
$A\times I\cup X\times \{0\}$ in $X\times I$
and a locally compact extension $G:V\to K$ of $H$.
Find a neighborhood $U$ of $A$ in $X$ such that $U\times I\subset V$
and pick a map $a:X\to I$ such that
$a(A)\subset \{1\}$ and $a(X-U)\subset \{0\}$.
Notice that $r:X\times I\to U\times I\cup X\times \{0\}$
defined by $r(x,t)=(x,t\cdot r(x))$ is continuous
and is identity on $A\times I\cup X\times \{0\}$.
Therefore the composition $H'=G\circ r$ is locally compact
and extends $H$.
\qed
\enddemo

Now we can reduce the question of extending a locally compact map
to the question of extending it up to homotopy to an arbitrary, not necessarily
locally compact, map.

\proclaim{2.13. Corollary} Suppose $X$ is a paracompact space,
$A$ is a closed subset of $X$, $K$ is
a CW complex, and $f:A\to K$ is a locally compact map. 
The following conditions are equivalent:
\item{a.} $f$ extends to a locally
compact map $f':X\to K$.
\item{b.} $f$ extends up to homotopy to a map $f':X\to K$.
\endproclaim
\demo{Proof} $a)$ is a special case of $b)$.
\par $b)\implies a)$. Suppose $f:A\to K$ is a locally compact map
and $g:X\to K$ is a map such that $g|A$ is homotopic to $f$.
Let $H:A\times I\cup X\times \{1\}\to K$ be a map such that
$H(x,0)=f(x)$ for $x\in A$ and $H(x,1)=g(x)$ for $x\in X$.
2.7 says that $H$ is homotopic to a locally compact map $H'$ in such a way
that the homotopy from $H$ to $H'$ is locally compact
on $A\times \{0\}$. Concatenating $H'$ with that homotopy
produces a locally compact $H'':A\times I\cup X\times \{1\}\to K$ such that
$H''(x,0)=f(x)$ for $x\in X$. By 2.12, $H''$ extends over $X\times I$
which gives a locally compact extension of $f$ over $X$.
\qed
\enddemo

\definition{2.14. Definition} $K$ is an absolute extensor up to homotopy
of $X$ if every map $f:A\to K$, $A$ closed in $X$, extends over $X$
up to homotopy.
\enddefinition

2.13 means that, if $X$ is paracompact and $K$ is a CW complex,
then $K\in AE_{lc}(X)$ is equivalent to $K$ being an absolute extensor of $X$
up to homotopy. Our next result relates the concept of being an absolute extensor
up to homotopy to the concept of being an absolute extensor in case
of simplicial complexes.

\proclaim{2.15. Theorem} Suppose $X$ is a paracompact space
and $K$ is a space.
Consider the following conditions:
\item{a.} $K$ is an absolute extensor of $X$
up to homotopy.
\item{b.} $K\in AE_0(X)$.
\item{c.} $K$ is an absolute extensor of $X$.
\par\noindent
If $K$ is an ANR for metrizable spaces, then Conditions a) and b)
are equivalent. If $K$ is complete ANR for metrizable spaces, then all three
conditions are equivalent.
\endproclaim
\demo{Proof} Assume $K$ is an ANR for metrizable spaces.
\par 
$a)\implies b)$. 
Suppose $f:A\to K$ is a map, where $A$ is a zero
subset of $X$.
Since $f$ extends over $X$ up to homotopy,
there is
 $H:A\times I\cup X\times \{1\}\to K$  such that
$H(x,0)=f(x)$ for $x\in A$.
Notice that $A\times I\cup X\times \{1\}$ is a zero subset of $X\times I$.
Therefore we can find a map $a:X\times I\to I$ such that
$A\times I\cup X\times \{1\}=a^{-1}(0)$.
Notice that $K$ can be considered as a subset of some Banach space
$E$. $E$ is an absolute
extensor of all paracompact spaces (see \cite{Hu}, Theorem 16.1b on p.63),
so there is an extension $G:X\times I\to E$ of $H$.
Consider the subset $K\times \{0\}\cup E\times (0,1]$
of $E\times I$.
Since $K$ is an absolute neighborhood extensor of all metrizable
spaces,
there is a retraction $r:U\to
K\times \{0\}$ from a neighborhood $U$ of $K\times \{0\}$
in $K\times \{0\}\cup E\times (0,1]$.
Define $F:X\times I\to K\times \{0\}\cup E\times (0,1]$
by $G'(x,t)=(F(x,t),a(x,t))$. $V=F^{-1}(U)$
is a neighborhood of $A\times I\cup X\times \{1\}$ is a closed subset of $X\times I$
and $r\circ F$ is an extension of $H$ over $V$.
Therefore $H$ extends over $X\times I$ which implies that $f$
extends over $X$.
\par $b)\implies a)$. Suppose $f:A\to K$ is a map
from a closed subset of $X$. Since $K$ is homotopy
equivalent to a CW complex, 2.6-2.7 and 2.11 imply
that there is a neighborhood $U$ of $A$ in $X$
and a homotopy extension
$f':U\to K$ of $f$.
Choose a map $a:X\to I$ such that $a(A)\subset \{0\}$ and
$a(X-U)\subset \{1\}$. Let $B=a^{-1}(0)$. $B$ is a zero subset
of $X$. Since $B\subset U$,
$f'|B$ extends over $X$ which proves that $K$ is an
absolute extensor of $X$ up to homotopy.
\par Assume $K$ is a complete ANR for metrizable spaces.
Obviously, Condition c) is stronger than Condition b).
\par $b)\implies c)$. Consider $K$ as a subset of
a Banach space $E$. Suppose $f:A\to K$ is a map
from a closed subset of $X$. Since $E$ is 
an absolute extensor of $X$,
there is an extension
$F:X\to E$ of $f$. Since $K$ is a $G_\delta$ subset of $E$,
$F^{-1}(K)$ is a $G_\delta$ subset of $X$ containing $A$.
Therefore there is a zero subset $B$ of $X$ so that
$A\subset B\subset F^{-1}(K)$. Now,
$F|B$ extends over $X$ which proves that $K$ is an
absolute extensor of $X$.
\qed
\enddemo

\head 3. Extension dimension for paracompact spaces \endhead

The purpose of this section is to prove existence of extension
dimension for paracompact spaces. It follows the same line of
reasoning as in \cite{Dr} for compact spaces or in \cite{D-D$_1$}
for separable metrizable spaces. The difference is that 2.9
allows for a significant simplification of the argument.

\proclaim{3.1. Proposition} Suppose $X$ is a paracompact space and $\{K_s\}_{s\in S}$ is
a family of pointed CW complexes. 
If each $K_s$ is an absolute extensor of $X$ up to homotopy,
then the wedge $K=\bigvee\limits_{s\in S}K_s$ is an absolute extensor
of $X$ up to homotopy.
\endproclaim
\demo{Proof} Let $K_T=\bigvee\limits_{s\in T}K_s$ for every finite subset
$T$ of $S$. $K_T\in AE_{lc}(X)$ for all $T$ implies $K\in AE_{lc}(X)$ by 2.9.
\qed
\enddemo

\proclaim{3.2. Proposition} Suppose $X$ is a paracompact space and $K\in AE_{lc}(X)$
is a CW complex. Let $n$ be the density of $X$ and let $m$ be a cardinal number
greater than or equal to $\max(2^n,2^{\aleph_0})$.
For any subcomplex $L$ of $K$ containing at most $m$ cells there is a subcomplex
$L'$ containing $L$ such that
\item{a.} $L'$ contains at most $m$ cells,
\item{b.} Any locally compact map $f:A\to L$, $A$ closed in $X$,
has a locally compact extension $f':X\to L'$.
\endproclaim
\demo{Proof} Let $Y$ be a dense subset of $X$ with cardinality equal to $n$.
Pick a point $\infty$ not belonging to $K$.
List all functions from $Y$ to $L\cup \{\infty\}$. There are at most
$m^n=m$ such functions. 
Keep only those functions $g$ so that
for some open set $U_g$ there is a locally
compact $u_g:cl(U_g)\to L$
so that $g(x)=u_g(x)$ for $x\in cl(U_g)\cap Y$
and $g(x)=\infty$ for $x\in Y-cl(U_g)$. Pick an extension $h_g:X\to K$
of $u_g$. The image $h_g(X)$ contains at most $m$ cells, so by adding all of them
we create a subcomplex $L'$ of $L$ containing at most $m$ cells.
\par
Any locally compact $f:A\to L$ extends over an open neighborhood
$U$ of $A$ in $X$. Let $f_1:U\to L$ be such extension
which is locally compact. Pick a neighborhood $V$ of $A$ in $X$
whose closure is contained in $U$.
Let $g:Y\to L\cup \{\infty\}$ be defined by $g(x)=f_1(x)$ if $x\in Y\cap cl(V)$,
$g(x)=\infty$ if $x\in Y-cl(V)$.
The function $g$ has a locally compact map $h_g:X\to K$
and $cl(U_g)\cap Y$ must be equal to $cl(V)\cap Y$.
Therefore $cl(U_g)=cl(V)$ and $h_g|A=f$. Thus, $f$ extends to a locally compact
map from $X$ to $L'$.
\qed
\enddemo

\proclaim{3.3. Corollary} Suppose $X$ is a paracompact space and $K\in AE_{lc}(X)$
is a CW complex. Let $n$ be the density of $X$ and let $m$ be a cardinal number
greater than or equal to $\max(2^n,2^{\aleph_0})$.
For any subcomplex $L$ of $K$ containing at most $m$ cells there is a subcomplex
$L'$ containing $L$ such that $L'$ contains at most $m$ cells
and $L'\in AE_{lc}(X)$.
\endproclaim
\demo{Proof} Put $L_1=L$. Create, using 3.2,
an increasing sequence of subcomplexes $L_n$ such that 
\item{a.} $L_n$ contains at most $m$ cells,
\item{b.} Any locally compact map $f:A\to L_n$, $A$ closed in $X$,
has a locally compact extension $f':X\to L_{n+1}$.
\par Apply 2.9 to the family $\{L_n\}_{n\ge 1}$
and conclude that $L'=\bigcup\limits_{n=1}^\infty L_n$
has the desired properties.
\qed
\enddemo

\proclaim{3.4. Proof of 1.3}
\endproclaim
 Let $n$ be the density of $X$ and let $m$ be the cardinal number
 equal to $\max(2^n,2^{\aleph_0})$.
Pick a set of CW complexes containing at most $m$ cells so that any
CW complex containing at most $m$ cells is listed there up to 
homeomorphism. 
Eliminate from that set CW complexes which are not absolute extensors of $X$
up to homotopy.
Let $\{K_s\}_{s\in S}$ be the resulting set
and put $K=\bigvee\limits_{s\in S}K_s$. By 3.1 $K$ is an absolute
extensor of $X$ up to homotopy.
Suppose $L$ is a CW complex which is an absolute extensor of $X$ up to homotopy.
We can express $L$ as the union of $\{L_t\}_{t\in T}$ of a partially ordered
family of subcomplexes of $L$ such that each $L_t$ is homeomorphic to one of $K_s$
(see 3.3).
If $K\in AE_{lc}(Y)$, then $K_s\in AE_{lc}(Y)$ for each $s\in S$ which implies
$L\in AE_{lc}(Y)$ by 2.9.
\qed

In practice one likes to be able to deal with absolute extensors rather than
absolute extensors up to homotopy. We are able to produce the extension dimension
of paracompact spaces by replacing CW complexes by complete simplicial 
complexes with the metric topology.
\proclaim{3.5. Proposition} For every CW complex $K$ there is
a simplicial complex $L$ such that $|L|_m$ is complete,
is homotopy equivalent to $K$,
and the following two conditions are equivalent for any paracompact space $X$:
\item{a.} $K$ is an absolute extensor of $X$ up to homotopy.
\item{b.} $|L|_m\in AE(X)$.
\endproclaim
\demo{Proof} Find a simplicial complex $M$ such that $|M|_m$
is homotopy equivalent to $K$ (see \cite{M-S}).
Triangulate $\bigcup\limits_{n=1}^\infty |M^{(n)}|_m\times [n,\infty)$
as $|L|_m|$ for some simplicial complex $L$.
Clearly, $|L|_m$ is homotopy equivalent to $K$. Suppose
it is an absolute extensor of $X$ up to homotopy.
Notice that $L$ does not contain any full infinite subcomplex.
Therefore $|L|_m$ is complete and 2.15
 implies that $|L|_m$ is an absolute extensor of $X$.
\qed
\enddemo

\proclaim{3.6. Proofs of 1.4 and 1.5}
\endproclaim
 By 1.3 there is a CW complex $K'$ such that
\item{1.} $K'$ is an absolute extensor of $X$ up to homotopy,
\item{2.} If a CW complex $L$ is an absolute extensor of
$X$ up to homotopy, then $L$ is an absolute extensor of
$Y$ up to homotopy of any paracompact space $Y$ such that
$K'$ is an absolute extensor of
$Y$ up to homotopy.
\par Pick a simplicial complex $K$ such that $|K|_m$ is complete,
is of the same homotopy type as $K'$, and $|K|_m\in AE(X)$ (see 3.5). 
\par 
Suppose $L$ is a complete ANR such that $L\in AE(X)$.
Choose a CW complex $L'$ of the same homotopy type as $L$.
Suppose $Y$ is a paracompact space such that $|K|_m\in AE(Y)$.
Now $K'$ is an absolute extensor of $Y$ up to homotopy
and $L'$ is an absolute extensor of $X$ up to homotopy.
Therefore $L'$ is an absolute extensor of $Y$ up to homotopy.
Since $L$ is homotopy equivalent to $L'$,
$L$ is an absolute extensor of $Y$ up to homotopy.
By 2.15, $L\in AE(Y)$.
\par 
Suppose $L$ is an ANR such that $L\in AE_0(X)$.
By 2.15, $L$ is an absolute extensor of $X$ up to homotopy.
Choose a CW complex $L'$ of the same homotopy type as $L$.
Suppose $Y$ is a paracompact space such that $|K|_m\in AE(Y)$.
Now $K'$ is an absolute extensor of $Y$ up to homotopy
and $L'$ is an absolute extensor of $X$ up to homotopy.
Therefore $L'$ is an absolute extensor of $Y$ up to homotopy.
Since $L$ is homotopy equivalent to $L'$,
$L$ is an absolute extensor of $Y$ up to homotopy.
By 2.15, $L\in AE_0(Y)$.

\qed

The Duality Theorem of Dranishnikov \cite{Dr} says
that each CW complex is equal to the extension dimension
of some compact Hausdorff space in the sense of Definition 1.1.
It is natural to ask if the same is true in the category
of paracompact spaces.
\proclaim{3.7. Problem} Suppose $K$ is a CW complex.
Is there a paracompact space $X$ so that $K$ is the extension dimension of $X$?
\endproclaim

An obvious approach to solve 3.7 is to produce a compact space
for $K$ as in \cite{Dr}. The remainder of this section is devoted
to explaining why this approach fails by showing paracompact spaces
whose extension dimension is not the same as of a compact space.

\definition{3.8. Definition} If $K$ and $L$ are
CW complexes, then $K\leq L$ means $L$ is an absolute
extensor up to homotopy of any paracompact space $X$
such that $K$ is an absolute extensor of $X$.
This leads to an equivalence relation $\sim$ on the category
of all CW complexes.
\par For any paracompact space $X$, $\ExD(X)$ stands
for its extension dimension in the sense of 1.3
and is unique up to equivalence $\sim$.
Now, for any paracompact spaces $X$ and $Y$, $X\leq Y$
means $\ExD(X)\leq \ExD(Y)$ and introduces a partial order
on the class of all paracompact spaces.
\enddefinition

Let us present a view of the Stone-\v Cech compactification
from the point of absolute extensors.

\proclaim{3.9. Proposition} 
In the class of normal spaces let $X\leq_f Y$ mean
that any finite CW complex $K$ which is an absolute
extensor of $Y$ must also be an absolute extensor of $X$.
Suppose $X$ is a normal space.
The class $\{Y\mid Y\leq_f X\text{ and Y is compact}\}$
has $\beta(X)$ as its maximum. Moreover, $X\leq_f\beta(X)$.
\endproclaim
\demo{Proof} 3.9 is well-known in the form:
$X$ and $\beta(X)$ have the same compact absolute extensors.
Let us sketch a proof for the sake of completeness.
Suppose $K\in AE(\beta(X))$. Any map $f:A\to K$, $A$ closed in $X$
extends over $\beta(A)$ which is a closed subset of $\beta(X)$.
Therefore $f$ extends over $\beta(X)$ and $K\in AE(X)$.
Suppose $K\in AE(X)$ and $f:A\to K$ is a map, $A$ closed in $\beta(X)$.
We can extend $f$ over a closed neighborhood $B$ of $A$
in $\beta(X)$. Let $g:X\to K$ be an extension of $f|B\cap X$.
Since $K$ is compact, $g$ extends over $\beta(X)$.
Let $h:\beta(X)\to K$ be such extension.
As $h$ and $g$ coincide on $\int(B)\cap X$,
they must coincide on $\int(B)$. In particular, $h$ is an extension of $f$.
\qed
\enddemo
Here is an extension theory analog of the Stone-\v Cech compactification. 
 
\proclaim{3.10. Theorem} Suppose $X$ is a paracompact space.
The class $\{Y\mid Y\leq X\text{ and Y is compact}\}$
has a maximum $X'$. There are separable metrizable
spaces $X$ such that
$X'<X$.
\endproclaim
\demo{Proof} Let $K$ be the extension dimension of $X$.
Let $X'$ be a compact Hausdorff space such that
$K\in AE(X')$ and $L\in AE(X')$, $L$ a CW complex,
implies $L\in AE(Y)$ for any compact Hausdorff space
$Y$ such that $K\in AE(Y)$ (see \cite{Dr}).
Since $K\in AE(X')$, $X'\leq X$. If $Y\leq X$
for some compact Hausdorff space $Y$,
then it simply means $K\in AE(Y)$.
To prove $Y\leq X'$ consider $M=\ExD(X')$. We need $M\in AE(Y)$
which follows from the way $X'$ was chosen.
\par 
In \cite{D-D$_2$}, Theorem 4.7, it is shown that if $G$ is a countable abelian group, 
and $A_p$ is the ring of $p$-adic integers for some prime number
$p$, then there is a separable space
$X$ of dimension 2 such that 
$\dim_GX\ne \dim_{A_p}\! X.$ Consider $G$ to be $Z$ localized at $p$
(all rational numbers with denominators relatively prime to $p$).
Now, $\ExD(X')=\ExD(X)$ implies
$\dim_GX'\ne \dim_{A_p}\! X'$ which is impossible for compact spaces
(see \cite{Ku}).
\qed
\enddemo

\head 4. Union theorem for paracompact spaces \endhead

In this section we prove the Union Theorem for paracompact spaces,
thus demonstrating that our extension theory of paracompact spaces
is quite natural.
\par To make sure that the approach in \cite{Dy$_1$} works we need the following result.

\proclaim{4.1. Lemma} Suppose $A$ is a subset of a hereditarily paracompact
space $X$. Any map $f:A\to K$ from $A$ to a CW complex $K$
extends up to homotopy
over a neighborhood of $A$ in $X$.
\endproclaim
\demo{Proof} It suffices to consider the case of $f$ being locally compact
and $K=|L|_w$ for some simplicial complex $L$.
Let $\{U_s\}_{s\in S}$ be a family of open sets in $X$ 
such that $A\subset U=\bigcup\limits_{s\in S} U_s$ and $f(A\cap U_s)$
is contained in a compact subset of $K$ for each $s\in S$.
Pick a locally finite partition $\{g_s\}_{s\in S}$ on $U$
($U$ is a paracompact space) such that $g_s^{-1}(0,1]\subset U_s$
for each $s\in S$. $\{g_s\}_{s\in S}$ may be viewed as a locally compact
map $g:U\to |L'|_w$, where $L'$ is the full simplicial complex
with the same vertices as $L$. Notice that $g|A$ is homotopic to $f$
as maps to $|L|_w$.
Pick a locally compact map $h:|L|_w\to |L|_w$ homotopic
to identity and extend it over a neighborhood $V$ of $|L|_w$
in $|L'|_w$.
Now, the composition of $g^{-1}(V)\to V\to |L|_w$
extends $f$ up to homotopy. \qed
\enddemo

\proclaim{4.2. Lemma} Suppose $A$ is an
$F_\sigma$-subset of a paracompact
space $X$. If $K$ is a CW complex which is an absolute extensor of $X$
up to homotopy, then $K$ is an absolute extensor of $A$ up to homotopy.
\endproclaim
\demo{Proof} $A$ is paracompact by 5.1.28 of \cite{En}.
Suppose $A=\bigcup\limits_{i=1}^\infty B_n$, where $B_n$ is a closed subset
of $X$ for each $n$. We may assume that $B_n\subset B_{n+1}$ for each $n$.
Suppose $C$ is a closed subset of $A$. Pick a closed subset
$D$ of $X$ such that $C=D\cap A$. Suppose $f:C\to K$ is a locally compact
map to a CW complex. 
Extend $f$ over a closed neighborhood $C_1$ of $C$ in $A$,
then use the fact that $K\in AE_{lc}(B_1)$ to extend it over
$C_1\cup B_1$. The resulting map $f_1:C_1\cup B_1\to K$
is locally compact by 2.4.
Suppose we have a locally compact map $f_n:C_n\cup B_n$
such that $C_n$ is a closed neighborhood of $C_{n-1}\cup B_{n-1}$ in $A$.
Extend it over a closed neighborhoof $C_{n+1}$ of $C_n\cup B_n$
and use the fact that $K\in AE_{lc}(B_{n+1})$ to extend it over
$C_{n+1}\cup B_{n+1}$. The resulting 
map $f_{n+1}:C_{n+1}\cup B_{n+1}\to K$ is locally compact by 2.4.
The direct limit $f'$ of maps $f_n$ is an extension of $f$ and is locally compact.
Indeed, given $x\in A$ we find the smallest $n$ such that $x\in C_{n}\cup B_{n}$.
$f'(x)$ equals $f_n(x)$. Since $C_{n+1}$ is a closed neighborhoof of $C_n\cup B_n$
and $f_{n+1}$ is locally compact, there is a neighborhood $U$ of $x$ in $A$
such that $f_{n+1}(U)=f'(U)$ is contained in a compact subset of $K$.
\qed
\enddemo

\proclaim{4.3. Theorem} Suppose $X$ is a hereditarily paracompact
space. Let $K$ and $L$ be CW complexes. If $K$ is an absolute
extensor of $A\subset X$ up to homotopy and $L$
is an absolute extensor of $B\subset X$ up to homotopy,
then the join $K\ast L$ is an absolute extensor of
$A\cup B$ up to homotopy. 
\endproclaim
\demo{Proof} It suffices to consider $X=A\cup B$.
We may assume that both $K$ and $L$ are
 simplicial complexes equipped with CW topology,
$K=|K'|_w$ and $L=|L'|_w$.
We will be working with locally compact maps
which are ideal for the following reason:
if $f:Y\to |M|_m$ is a map such that
every $y\in Y$ has a neighborhood $U$ with $f(U)$
contained in a finite subcomplex of $|M|_m$,
then $f$ considered as a function from $Y$ to $|M|_w$
is continuous.
\par

 Suppose $C$ is a closed subset of $A\cup B$
and $f:C\to K*L$ is a locally compact map.  Notice that $f$ defines two closed,
disjoint subsets $C_K=f\1 (K)$,  $C_L=f\1 (L)$ of
$C$ and locally compact
 maps $f_K:C-C_L\to K$, $f_L:C-C_K\to L$, $\alpha :C\to
[0,1]$ such that: 
\item{1.} $\alpha \1 (0)=C_K$, $\alpha \1
(1)=C_L$, 
\item{2.} $f(x)=(1-\alpha (x))\cdot f_K(x)+\alpha (x)\cdot f_L(x)$ for
all $x\in C$. 
\par\noindent
Indeed, each point $x$ of a simplicial complex $M$ can be uniquely
written as $x=\sum_{v\in M^{(0)}}\phi_v(x)\cdot v$,
where $M^{(0)}$ is the set of vertices of $M$ ($\{\phi_v(x)\}$ are called
{\it barycentric coordinates} of $x$). We define
$\alpha (x)$ as $\sum_{v\in L^{(0)}}\phi_v(f(x))$,
$f_K(x)$ is defined as $(\sum_{v\in K^{(0)}}\phi_v(f(x))\cdot  v)/(1-\alpha (x))$
and $f_L(x)$ is defined as $(\sum_{v\in L^{(0)}}\phi_v(f(x))\cdot  v)/(\alpha (x))$. 
\par Since $K\in AE_{lc}(A-C_L)$ by  4.2, 
$f_K$ extends over $(C\cup A)-C_L$. To make sure that there is a locally compact extension
we proceed as follows: first extend $f_K$ over a closed neighborhood
$D$ of $C-C_L$ in $(C\cup A)-C_L$. Let $u:B\to K$
be a locally compact extension  of $f_K|(C-C_L)$.
Extend $u|B\cap (A-C_L)$ to a locally compact $v:A-C_L\to K$.
Pasting $v$ and $f_K$ results in a locally compact map.
\par

 Consider a homotopy extension  $g_K:U_A\to K$ of $f_K$ over a
neighborhood $U_A$ of $(C\cup A)-C_L$ in $X-C_L$. Since
$C-C_L$ is closed in $U_A$, we may assume that $g_K$ is an actual
extension of $f_K:C-C_L\to K$ (see 2.13).
Similarly, let  $g_L:U_B\to L$ be
an extension of $f_L$ over a neighborhood $U_B$ of $(C\cup B)-C_K$ 
in $X-C_K$. Notice that $X=U_A\cup U_B$. Let $\beta
:X\to [0,1]$ be an extension of $\alpha $ such that
$\beta(X-U_B)\subset \{0\}$ and $\beta(X-U_A)\subset \{1\}$. 
 Define  $f':X\to K*L$ by
$$f'(x)=(1-\beta (x))\cdot g_K(x)+\beta (x)\cdot g_L(x)\text{ for all } 
x\in U_A\cap U_B,$$
$$f'(x)=g_K(x)\text{ for all } 
x\in U_A-U_B,$$
and $$f'(x)=g_L(x)\text{ for all } 
x\in U_B-U_A.$$ 
Notice that $f'$ is an
extension of $f$.  
Now, it suffices to prove that $f':X\to |K'\ast L'|_m$
is continuous. Indeed, as identity $|K'\ast L'|_w\to |K'\ast L'|_m$
is a homotopy equivalence it would certify the existence
of an extension of $f:C\to |K'\ast L'|_w$ up to homotopy
which is all we need in view of 2.13.
\par
To prove the continuity of $f':X\to |K'\ast L'|_m$ we need to show that $\phi_vf'$ is continuous
for all vertices $v$ of $K'*L'$ (see \cite {M-S, Theorem 8 on p.301}).
Without loss of generality we may assume that $v\in K'$. Then, 
$$\phi_vf'(x)=(1-\beta (x))\cdot \phi_vg_K(x) \text{ for all } 
x\in U_A$$
and $$\phi_vf'(x)=0\text{ for all } 
x\in U_B-U_A.$$ 
Clearly, $\phi_vf'|U_A$ is continuous. Suppose $x_0\in (U_B-U_A)\cap cl(U_A)$
and $M>0$. Since $\phi_vf'(x_0)=0$, it suffices to show existence
of a neighborhood $W$ of $x_0$ such that $\phi_vf'(W)\subset [0,M)$.
As $\beta(x_0)=1$, there is a neighborhood $W$ of $x_0$
so that $\beta(W)\subset (1-M,1]$.
If $x\in W\cap (U_B-U_A)$, then $\phi_vf'(x)=0$.
If $x\in W\cap U_A$, then $\phi_vf'(x)=(1-\beta(x))\cdot \phi_vg_K(x)\leq
1-\beta(x)<M$.
\qed
\enddemo

\head 5. Spaces with all maps being locally compact \endhead

It is of interest to see which maps to CW complexes are locally compact.

\proclaim{5.1. Problem} Characterize all paracompact spaces $X$ so that 
any map $f:A\to K$, $A$ closed in $X$ and $K$ a CW complex, is locally
compact.
\endproclaim
This section is devoted to partial answers to 5.1.

\proclaim{5.2. Proposition} Suppose $f:X\to Y$ is a perfect map
and $X$ is a paracompact space.
If every map from $X$ to a CW complex is locally compact,
then every map from $Y$ to a CW complex is locally compact.
\endproclaim
\demo{Proof} Suppose $g:Y\to K$ is a map from $Y$ to a CW complex.
Let $y_0\in Y$. Since $g\circ f$ is locally compact,
for each $x\in f^{-1}(y_0)$ there is a neighborhood $U_x$
such that $gf(U_x)$ is contained in a compact subset $Z_x$ of $K$.
As $f^{-1}(y_0)$ is compact, $f^{-1}(y_0)\subset \bigcup\limits_{x\in F}U_x$
for some finite subset $F$ of $f^{-1}(y_0)$.
Since $f$ is closed there is a neighborhood $U$ of $y_0$ in $Y$
with $f^{-1}(U)\subset \bigcup\limits_{x\in F}U_x$.
Now $g(U)=gf(f^{-1}(U))\subset gf(\bigcup\limits_{x\in F}U_x)\subset \bigcup\limits_{x\in F}Z_x$
which proves that $g$ is locally compact.
\qed
\enddemo

\proclaim{5.3. Proposition} Suppose $A$ is a subset of $X$
and has a countable basis of neighborhoods. If $f:X\to K$
is a map to a CW complex such that $f(A)$ is contained
in a compact subset of $K$, then there is a neighborhood $U$
of $A$ in $X$ such that $f(U)$ is contained in a compact subset of $K$.
\endproclaim
\demo{Proof} There is a finite subcomplex $K_0$ of $K$ containing $f(A)$.
Choose a basis of neighborhoods $\{U_n\}_{n\ge 1}$ of $A$ in $X$.
Suppose none of $f(U_n)$ is contained in a finite subcomplex of $K$.
Choose, by induction, elements $w_n\in f(U_n)$ so that
the smallest subcomplex of $K$ containing $K_0$ and $w_1,\ldots,w_{n-1}$
does not contain $w_n$.
The set $C=\{w_i\}_{i\ge 1}$ is closed in $K$ and misses $K_0$,
so $f^{-1}(C)$ is closed and misses $A$. Pick $m$ so that $U_m\subset X-f^{-1}(C)$.
Now $w_m\in K-C$,
a contradiction.
\qed
\enddemo

\proclaim{5.4. Corollary} If $X$ is the union of its compact subsets
which have a countable basis of neighborhoods, then any map from $X$
to a CW complex is locally compact.
\endproclaim

\remark{Remark} Hausdorff spaces $X$ such that every point is contained
in a compact subset $Z$ with countable basis of neighborhoods
are discussed in \cite{En} (Exercise 3.1.E to section 1 of chapter 3)
under the name of {\bf pointwise countable type}.
The class of such spaces contains locally compact spaces,
first countable spaces, is closed under finite cartesian products,
is hereditary with respect to closed
subsets, and is hereditary with respect to G$_\delta$-subsets
(in particular, all topologically complete spaces belong to the class).
It is also easy to show that if $f:X\to Y$ is a perfect map
and $Y$ belongs to the class, than $X$ belongs to the class.
\endremark

\proclaim{5.5. Problem} Suppose $X$ is a paracompact space such that any map
from a closed subset $A$ of $X$ to a CW complex is locally compact.
Let be $Y$ a compact
space. Is every map from a closed subset $A$ of $X\times Y$
to a CW complex locally compact?
\endproclaim

\enddocument